\documentclass[reqno]{amsart}
\usepackage{amsthm,amsmath,amssymb}
\usepackage[pdftex]{graphicx}
\usepackage{braket}
\usepackage{color}
\usepackage{cleveref}
\usepackage{amscd}
\numberwithin{equation}{section}

\title{Counterexamples to Allen's conjectures}
\author{Kouki Sato}
\date{}
\address{\textsc{Meijo University,Tempaku, Nagoya 468-8502, Japan}}
\email{satokou@meijo-u.ac.jp}

\newtheorem{thm}{Theorem}%[section]

\newtheorem{lem}[thm]{Lemma}
\newtheorem{cor}[thm]{Corollary}

\newtheorem{question}[thm]{Question}
\newtheorem{conj}[thm]{Conjecture}

\theoremstyle{definition}

\newtheorem{remark}[thm]{Remark}
\newtheorem*{acknowledge}{Acknowledgements}

\begin{document}

\begin{abstract}
We show that the torus knots $T(2,5)$ and $T(2,9)$ bound smooth M\"{o}bius bands
in the 4-ball whose double branched covers are negative definite,
giving counterexamples to Conjectures 1.6 and 1.8 of Allen in [New York J.\ Math.\ 29 (2023) 1038--1059]. 
\end{abstract}

\maketitle

\section{Introduction}

In \cite{Allen:2023}, Allen asked the following question. 

\begin{question}[\text{\cite[Question 1.2]{Allen:2023}}]
\label{ques:allen}
Given a knot $K \subset S^3$, what is the set of realizable pairs
\[
(e(F),h(F)) = (\text{normal Euler number of a surface $F$}, \text{first Betti number of $F$}),
\] 
where $F \subset B^4$ is bounded by $K$?
\end{question}

She studied \Cref{ques:allen} for the torus knots $T(2,n)$ and $T(3,n)$ with $n>0$,
and proved the following result about $T(2,n)$.

\begin{thm}[\text{\cite[Theorem 1.4]{Allen:2023}}]
\label{thm:allen}
For $T(2,n)$ with odd $n>0$, the following pairs are realizable:
\[
(e,h) \in \{(-2n \pm 2m, 1+m+2l) \mid m,l \geq 0\} \cup \{(2+2m, n+m) \mid m \geq 0\},
\]
it is unknown if the following pairs are realizable:
\begin{enumerate}
\item if $n \equiv 1$ (\text{mod $4$}), $(e,h) = (4-2n + 2m, 1+m)$ for $0 \leq m < n-1$, 
\item if $n \equiv 3$ (\text{mod $4$}), $(e,h) = (8-2n + 2m, 3+m)$ for $0 \leq m < n-3$, 
\end{enumerate}
and all other pairs are not realizable.
\end{thm}
\begin{remark}
In the original statement of \cite[Theorem 1.4]{Allen:2023},
the pair $(2,n)$ is included in both the set of realizable pairs and the set of unknown points, but actually the pair is realized by the boundary connected sum of the minimal genus Seifert surface for $T(2,n)$ and the standard
M\"{o}bius band $M$ with $e(M)=2$. 
\end{remark}

For the unknown points in \Cref{thm:allen}, she posted the following conjecture.

\begin{conj}[\text{\cite[Conjecture 1.6]{Allen:2023}}]
\label{conj1}
All unknown points in \Cref{thm:allen} are not realizable.
\end{conj}

For general torus knots, the following conjecture was also posted.

\begin{conj}[\text{\cite[Conjecture 1.8]{Allen:2023}}]
\label{conj2}
All torus knots have a single realizable “minimal point". In other words, for a torus knot $K$, there is exactly one realizable pair of the form $(e, \gamma_4(K))$.
\end{conj}
Here $\gamma_4(K)$ is the {\it non-orientable 4-genus} of $K$, i.e.\ the minimal first Betti number of non-orientable surfaces in $B^4$ with boundary $K$.

It follows from the Gordon-Litherland formula $\sigma(\Sigma(F)) = \sigma(K) - e(F)/2$ in \cite{GL:1978} that the unknown points in \Cref{thm:allen} are corresponding to surfaces whose double branched covers are negative definite and first Betti numbers are less than $2g_4(T(2,n))+1$. Here $\Sigma(F)$ is the double branched cover of $B^4$ over $F$, $\sigma(K)$ is the knot signature of $K$ and
$g_4(K)$ is the {\it 4-genus} of $K$, i.e.\ the minimal genus of orientable surfaces in $B^4$ with boundary $K$. (Precisely, the cases where $n \equiv 3 \mod 4$ and $h(F) \in \{1,2\}$ are ruled out by \Cref{thm:allen}.)
Also note that $T(2,n)$ obviously bounds a smooth M\"{o}bius band whose double branched cover is positive definite.
In this paper, we show the following.
\begin{thm}
\label{thm:main}
The torus knots $T(2,5)$ and $T(2,9)$ bound smooth M\"{o}bius bands
in the 4-ball whose double branched covers are negative definite.
In particular, these surfaces realize $(-6, 1)$ and $(-14,1)$ respectively, giving counterexamples to Conjectures \ref{conj1} and \ref{conj2}. 
\end{thm}
Since the unknot bounds a smooth M\"{o}bius band $M$ with $e(M)= 2$,
the realizability of $(e,h)$ implies the realizability of $(e + 2m, h+ m)$ for any $m \geq 0$.
Moreover, since there exists an orientable cobordism $C$ from $T(2,n)$ to $T(2,n+2)$ in $S^3\times[0,1]$ with $b_1(C)=2$, the realizability of $(e,h)$ for $T(2,n)$ implies the realizability of $(e, h+2)$ for $T(2,n+2)$.
Therefore, we have the following corollary.
\begin{cor}
\label{cor}
For $T(2,n)$ with $n=5,7,9,11$, all unknown points in \Cref{thm:allen} are realizable.
For $T(2,n)$ with $n \geq 13$, the following pairs are realizable:
\begin{enumerate}
\item if $n=4k+9$, $(e,h) = (4-2n + 2m, 1+m)$ for $4k \leq m < n-1$, and
\item if $n=4k+11$, $(e,h) = (8-2n + 2m, 3+m)$ for $4k \leq m < n-3$. 
\end{enumerate}
\end{cor}

In the light of \Cref{cor}, it would be natural to ask the following question.

\begin{question}
Are the all unknown points in \Cref{thm:allen} realizable?
\end{question}

\begin{acknowledge}
The author would like to thank Joshua Sabloff 
and Masaki Taniguchi for their stimulating conversations.
\end{acknowledge}

\section{Counterexamples}

\begin{figure}[tbp]
\includegraphics[scale=0.7]{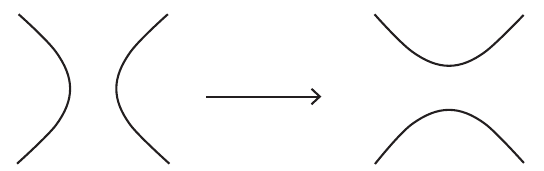}
\caption{\label{fig:H(2)} $H(2)$-move}
\end{figure}

The {\it $H(2)$-move} (or {\it pinch move}) is a deformation of link diagram shown in \Cref{fig:H(2)}.
If a diagram $D_1$ for a knot $K_1$ is deformed into a diagram $D_2$ for a knot $K_2$ by an $H(2)$-move,
then it gives rise to a non-orientable cobordism $C$ from $K_1$ to $K_2$ in $S^3 \times [0,1]$
such that $b_2(\Sigma(C))=1$. The computation of $\sigma(\Sigma(C))$ is given by the following formula. (Here $w(D_i)$ denotes the writhe of $D_i$.)
\begin{lem}[\text{\cite[Lemma 6.6]{Sato:2019}}]
\label{lem:signature}
$
\sigma(\Sigma(C))= \big(\sigma(K_2) -\sigma(K_1) \big) + \frac{1}{2} \big(w(D_2) -w(D_1) \big).
$
\end{lem}

\begin{figure}[tbp]
\includegraphics[scale=0.95]{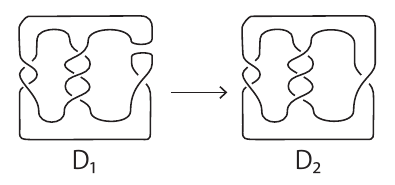}
\caption{\label{fig:T25} $H(2)$-move from the unknot to $T(2,5)$}
\end{figure}

\begin{figure}[tbp]
\includegraphics[scale=0.95]{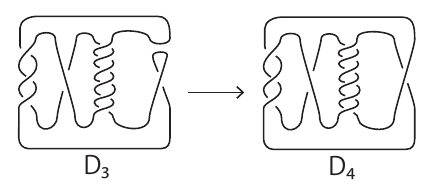}
\caption{\label{fig:T29} $H(2)$-move from $6_1$ to $T(2,9)$}
\end{figure}

\begin{proof}[Proof of \Cref{thm:main}]
Let $C \colon K_1 \to K_2$ be a cobordism in $S^3 \times [0,1]$ derived from an $H(2)$-move shown in \Cref{fig:T25}.
Then, it is easy to check that $D_1$ and $D_2$ are diagrams for the unknot $U$ and $T(2,5)$
respectively, and
\begin{align*}
\sigma(\Sigma(C)) &= \big(\sigma(K_2) -\sigma(K_1) \big) + \frac{1}{2} \big(w(D_2) -w(D_1) \big)\\[2mm]
&= -4 + \frac{1}{2}(6-0) = -1.
\end{align*}
Gluing $C$ with a standard disk in  $B^4$ along $K_1=U$, we have a smooth M\"{o}bius band $F$
with boundary $K_2=T(2,5)$ and $b_1(F)=b_2(\Sigma(F))=-\sigma(\Sigma(F))=1$.

Similarly, let $C' \colon K_3 \to K_4$ be a cobordism in $S^3 \times [0,1]$ derived from an $H(2)$-move shown in \Cref{fig:T29}.
Then, $D_3$ and $D_4$ are diagrams for the knot $6_1$ in the Rolfsen table \cite{Rolfsen:1972} and $T(2,9)$ respectively, and
\begin{align*}
\sigma(\Sigma(C')) &= \big(\sigma(K_4) -\sigma(K_3) \big) + \frac{1}{2} \big(w(D_4) -w(D_3) \big)\\[2mm]
&= -8 + \frac{1}{2}(9-(-5)) = -1.
\end{align*}
Gluing $C'$ with a slice disk for $6_1$ in  $B^4$ along $K_3=6_1$, we have a smooth M\"{o}bius band $F'$
with boundary $K_4=T(2,9)$ and $b_1(F')=b_2(\Sigma(F'))=-\sigma(\Sigma(F'))=1$.
\end{proof}


\bibliographystyle{plain}
\bibliography{tex}

\end{document}